\documentclass[reqno]{amsart}
\usepackage{etex}
\usepackage{graphicx}
\usepackage{amsmath}
\usepackage{mathrsfs}
\usepackage{amssymb}
\usepackage[all]{xy}
\usepackage{color}
\usepackage{pictex}
\usepackage{subfigure}
\usepackage{xcolor}
\usepackage[colorlinks]{hyperref}

\newtheorem{thm}{Theorem}[section]
\newtheorem{cor}[thm]{Corollary}
\newtheorem{lem}[thm]{Lemma}
\newtheorem{prop}[thm]{Proposition}
\theoremstyle{definition}

\theoremstyle{remark}

\newtheorem{rem}[thm]{Remark}
\numberwithin{equation}{section}

\newcommand{\z}{\mathbb Z}
\newcommand{\cp}{\mathbb C \mathrm P}
\newcommand{\hp}{\mathbb H \mathrm P}
\newcommand{\m}{\mathrm{MCG}}
\newcommand{\e}{\mathrm{haut}}
\newcommand{\T}{\mathcal T}
\newcommand{\omsp}{\Omega^{\mathrm{spin}}}

\newcommand{\om}{\Omega}
\newcommand{\ov}{\overline}
\newcommand{\wt}{\widetilde{\times}}
\newcommand{\dt}{\mathbf t}
\begin{document}

\title[homotopy automorphism group of homology $\mathbb C \mathrm P^3$]{On the homotopy automorphism group of simply-connected homology $\mathbb C \mathrm P^3$}
\author{Matthias Kreck, Yang Su, Shengye Zhu}
\address{Mathematisches Institut, Universit\"at Bonn}
\email{kreck@math.uni-bonn.de}
\address{HLM, School of Mathematics and Systems Science, Chinese Academy of Sciences; 
School of Mathematical Sciences, University of Chinese Academy of Sciences, Beijing 100190, China}
\email{suyang@math.ac.cn}
\address{School of Mathematical Sciences, University of Chinese Academy of Sciences, Beijing 100049, China}
\email{zsy97@amss.ac.cn}


\date{}
\begin{abstract}
In this paper, we compute the group of homotopy equivalences of closed, smooth, simply connected $6$-manifolds whose integral homology groups are isomorphic to whose of the complex projective space $\cp^3$, building on computations of their mapping class groups.
\end{abstract}
\maketitle
\section{Introduction}
For a topological space $X$, the group of homotopy equivalences $\e(X)$ is an object of significant interest in homotopy theory, and has been widely studied in the literature; see, for example, the survey articel \cite{Ar90}. For a compact manifold $M$, the group $\e(M)$ also plays an important role in manifold topology, particularly in the classification of manifolds and the study of group actions, see \cite{Ba96}. 

For simply connected $4$- and $5$-manifolds, the group $\e(M)$ has been computed by Milnor \cite{Mil} and Baues-Buth \cite{BB96}, respectively.  In dimension $6$, however, the situation is less well understood.  The group $\e(S^2 \times S^4)$ was computed by Sawashita (\cite{Sawa}) as part of a broader study of the group of homotopy equivalences of products of spheres. Partial information on the group of homotopy equivalences of $S^2$-bundles over $S^4$ has been obtained by Yamaguchi (\cite{Yama}).

As a generalization of the results mentioned above,  in this paper we compute the group $\e(M)$ for closed, smooth, simply connected $6$-manifolds $M$ with $H_*(M) \cong H_*(\cp^3)$ (referred to as simply connected homology $\cp^3$ for short). For this class of manifolds, classifications up to diffeomorphisms and up to homotopy equivalences have long been known. We briefly recall these results below. 

By the theory of Wall (\cite[Theorem 3]{Wall}) and Jupp (\cite[Theorem 1]{Jupp}), the diffemorphism type of an oriented, simply connected $6$-manifold $M$ with torsion-free homology is determined by the third Betti number $b_3(M)$, the cubic form 
$$H^2(M) \times H^2(M) \times H^2(M) \to \z, \ \ (x, y, z) \mapsto \langle x \cup y \cup z, [M] \rangle,$$ 
the second Stiefel-Whitney class   $w_2(M)$ and the first Pontrjagin class $p_1(M)$. 

When $M$ is a homology $\cp^{3}$, the cubic form and $p_1(M)$ can be described in terms of numeric invariants as follows. Let $x \in H^2(M)$, $y \in H^{4}(M)$ be generators. Then there exists an integer $d$ such $x^2 = d y$.  If $M$ is nonspin, then $d$ is even. By choosing the sign of $y$, we may assume $d \ge 0$. 
There exists an integer $l$ such that the first Pontrjagin class is given by
 $$p_1(M)=4(d+6l) y,  \ \textrm{if $M$ is spin;} \ \ \ p_1(M)=(d+48l)y,  \ \textrm{if $M$ is non-spin}.$$ 
 When $d=0$, we choose the sign of $y$ so that $l \ge 0$. 
 
 The homotopy classification of this class of $6$-manifolds has been established by several authors (see \cite[Proposition 1.9]{Baues} and the discussion therein).
We summarize these results in the theorem below.

\begin{thm}\label{thm:class}
Let $M$ be a simply connected homology $\cp^3$. Then $M$ admits an orientation-reversing diffeomorphism.

If $M$ is spin, there is a one-to-one correspondence between the set of diffeomorphism classes of such manifolds and the set 
$$\{ (d, l) \in \z \times \z \ | \ d \ge 0, \textrm{ and $l \ge 0$ when $d=0$} \} .$$ 
If $d$ is odd, the homotopy type of $M$ is determined by $d$; if $d$ is even,  the homotopy type of $M$ is determined by $d$ together with $d/2+l \pmod 2$. 

If $M$ is non-spin, there is a one-to-one correspondence between the set of diffeomorphism classes and the set 
$$\{ (d, l) \in 2\z \times \z \ | \ d \ge 0, \textrm{ and $l \ge 0$ when $d=0$} \} .$$ 
In this case, the homotopy type of $M$ is determined by $d$.

\end{thm} 

We now give some example. In the spin case, the manifold with $d=l=0$ is $S^2 \times S^4$, while the manifold with $d=0$ and $l=1$ is $\cp^3$; in the non-spin case, the manifold with $d=l=0$ is the nontrivial $S^4$-bundle over $S^2$, denoted by $S^4 \wt S^2$. 

Another class of geometrically interesting examples is given by the total spaces $X_n$ of $S^2$-bundle over $S^4$ with characteristic map $n \in \pi_3SO(3) \cong \z$. These manifolds are spin and satisfy $d=n$ and $l=0$.   As a consequence of the homotopy classification, if $M$ is spin with $d$ odd, or with both $d$ and $l$ even, then $M$ is homotopy equivalent to  the sphere bundle $X_d$.

\smallskip

A standard method in homotopy theory for studying the group of homotopy equivalences of an $n$-dimensional simply connected Poincar\'e complex $X$ is provided by the exact sequence of Barcus and Barratt \cite[Theorem 6.1]{BB58} (see also \cite[Proposition (1.3)]{Ba96}) 
$$\mathrm{Im}\, \{ \pi_n(X^{(n-1)}) \to \pi_n(X)\} \stackrel{d^*}{\longrightarrow}  \e(X) \stackrel{j^*}{\longrightarrow}  \e(X^{(n-1)}),
$$
where $X^{(n-1)}$ denotes the $(n-1)$-skeleton of $X$. The homomorphism 
$$d^* \colon \mathrm{Im}\, \{ \pi_n(X^{(n-1)}) \to \pi_n(X)\} \to \e(X)$$ 
is induced by the composition of the coaction $\mu \colon X \to X \vee S^n$ with a map $S^n \to X$ representing an element of $\mathrm{Im}\, \{ \pi_n(X^{(n-1)}) \to \pi_n(X)\} $. There are clear descriptions of $\mathrm{Im} \, j^*$ and the action of $\mathrm{Im} \, j^*$ on $\mathrm{Im} \, d^*$. The difficulty lies in the determination of the kernel of $d^*$, which involves the computation of a certain homomorphism $\pi_1(\mathrm{Map}_* (X^{(n-1)},X), i) \to \pi_n(X)$. 

This method has been applied to the class of $6$-manifolds under consideration by several authors. For $X= S^2 \times S^4$, Sawashita showed that this sequence extends to a split short exact sequence (\cite[Theorem 4.3]{Sawa})
$$0 \to \pi_6(S^2 \times S^4) \to \e(S^2 \times S^4) \to (\z/2)^3 \to 1.$$

For $X=X_d$, Yamaguchi obtained an exact sequence (\cite[Theorem 1.1]{Yama})
\begin{equation}\label{Yama}
\pi_6(X_d) \to \e(X_d) \to G_d \to 1
\end{equation}
with $G_d \cong \z/2$ when $d$ is odd, and $G_d \cong \z/2 \times \z/2$ when $d$ is even. He also computed the group $\pi_6(X_d)$: 
\begin{equation}\label{eq:pi6}
\pi_6(X_d) \cong \left \{ \begin{array}{ll}
0 & \mathrm{if} \ (d,6)=1 \\
\z/3 & \mathrm{if} \ (d,6)=3 \\
\z/2 \times \z/(12, d/2) & \textrm{if $d$ is even} 
\end{array} \right. 
\end{equation}
However, since the kernel of the homomorphism $\pi_6(X_d) \to \e(X_d)$ was not determined (except for the cases where $\pi_6(X_d)=0$), the problem of computing $\e(X_d)$ remained unsolved.

In this paper, we study the group of homotopy equivalences from a different perspective. In manifold topology, the mapping class group $\m(M)$ of a smooth manifold $M$ plays an analogous role to that of the group of homotopy equivalences in homotopy theory,  and is itself an object of considerable interest. The group $\e(M)$ is related to group $\m(M)$ via the natural forgetful homomorphism $\psi \colon \m(M) \to \e(M)$. 

As an example,  for $(n-1)$-connected $2n$-dimensional almost parallelizable manifolds, both $\m(M)$ and $\e(M)$ have been computed by various authors, and a detailed comparison between them is given in \cite[\S 10]{Ba96}). 

In this paper, we analyze the forgetful homomorphism $\psi$ for simply connected homology $\cp^3$. We determine the kernel and cokernel of $\psi \colon \m(M) \to \e(M)$, and then use our knowledge of $\m(M)$, obtained in \cite{KrSu}, to compute the group $\e(M)$. The manifold-theoretic perspective adopted in this paper also allows us to describe the generators of \(\e(M)\) geometrically. 

Before stating the main theorems, we present the construction of some self-homotopy equivalences of $M$, which serve as generators of $\e(M)$. Let $\mu \colon M \to M \vee S^6$ be the pinch map which collapses the boundary of an embedded $6$-disk to a point. Let $\xi \colon S^6 \to M$ be the composite of maps $S^6 \to S^2$ and $S^2 \to M$ representing a generator of $\pi_6(S^2) \cong \z/12$ and $\pi_2(M) \cong \z$, respectively.  
When $d$ is even, it is shown in \cite[Proposition 5.3]{Yama} that $\pi_6(X_d)$ contains a direct summand which is mapped isomorphically to $\pi_6(S^4) \cong \z/2$ under the bundle projection $X_d \to S^4$. 
When $M$ is spin with both $d$ and $l$ even, by the homotopy classification, $M$ is homotopy equivalent to $X_d$. 
Let $\eta \colon S^6 \to M$ represent the corresponding non-trivial element in $\pi_6(M)$. 
Define homotopy equivalences $f_{\xi}$ and $f_{\eta}$
$$f_{\xi} = (\mathrm{id} \vee \xi) \circ \mu \colon  M \to M \vee S^6 \to M, \ f_{\eta} = (\mathrm{id} \vee \eta) \circ \mu \colon  M \to M \vee S^6 \to M.$$
In \S \ref{subsec:dehn}, we define Dehn twists
$$\dt_{\epsilon}, \  \dt_{\rho} \colon M \to M,$$ 
which are diffeomorphisms of $M$ supported in a tubular neighborhood of an embedded $S^3$ in $M$. 

\begin{prop}\label{prop}
$f_{\xi}$  is homotopic to  \(\dt_{\epsilon}\); $f_{\eta}$ is not homotopic to a diffeomorphism. 
\end{prop}
For integers $m$ and $n$, we write $(m,n)$ for their greatest common divisor. 

\begin{thm}\label{thm:mainspin}
Let $M$ be a simply connected spin homology $\cp^3$ with $d \ne 0$. 

If $d$ is odd, then 
$$\e(M) \cong \z/(3,d) \rtimes \z/2,$$
the generator of $\z/(3,d)$ is represented by $f_{\xi}$, the generator of $\z/2$ is represented by an orientation-reversing self-diffeomorphism of $M$, and the action of $\z/2$ on $\z/(3,d)$ is by the multiplication of $-1$.
 
If $d$ is even,  then there is a short exact sequence
 \begin{equation}\label{eq:deven}
 1 \to G \rtimes \z/2  \to \e(M) \to \mathbb Z/2  \to 1,
 \end{equation}
where 
\[
G \cong
\begin{cases}
\mathbb Z/2 \times \mathbb Z/(6,d/2)
& \text{if } d \not\equiv 0 \pmod 8 \text{ and } l \not\equiv 0 \pmod 2, \\[0.3em]
\mathbb Z/2 \times \mathbb Z/(12,d/2)
& \text{otherwise}.
\end{cases}
\]
In the group $G$, the generator of $\z/2$  is represented by $\dt_{\rho}$,  the generator of the other direct summand is represented by $f_{\xi}$. The generator of factor $\z/2$ in the semi-direct product is represented by an orientation-reversing smooth conjugation on $M$. The action of $\z/2$ on $G$ is by the multiplication by $-1$.
If $l$ is even, then the exact sequence (\ref{eq:deven}) splits, with a splitting given by $f_{\eta}$.
\end{thm}

\begin{rem}
Together with Sawashita's computation of $\e(S^2 \times S^4)$, the above theorem provides a description of $\e(M)$ for all simply connected spin homology $\cp^3$'s, except for the unique homotopy type with $d=0$ and $l$ odd (see Theorem \ref{thm:class}). In the remaining case, we show that $\e(M)$ contains a subgroup of index $2$ or $4$ isomorphic to $(\z/2 \times \z/6)\rtimes \z/2$; see Theorem \ref{thm:hmtp1} (5). 
\end{rem}

As a corollary, the exact sequence (\ref{Yama}) extends to a split short exact sequence. 
\begin{cor}\label{cor}
Let $X_d$ be the total space of the $S^2$-bundle over $S^4$ with characteristic class $d \in \pi_3SO(3)$. Assume that $d \neq 0$. Then there is a split short exact sequence
$$0 \to \pi_6(X_d) \to \e(X_d) \to G_d \to 1.$$  
\end{cor}

From this perspective, the generators of the corresponding groups are described as follows:
when \(d\) is odd, a generator of \(\pi_6(X_d)\) is represented by the Dehn twist \(\dt_{\epsilon}\), and a splitting of \(G_d \cong \mathbb Z/2\) is given by the fiberwise antipodal map \(f_a\). When \(d\) is even, a generator of the direct summand $
\mathbb Z/(12,d/2) < \pi_6(X_d)$
is represented by the Dehn twist \(\dt_{\epsilon}\), while a generator of the direct summand
$
\mathbb Z/2 < \pi_6(X_d)
$
is given by $f_{\eta}$. A splitting of
$
G_d \cong \mathbb Z/2 \times \mathbb Z/2
$
is provided by \(f_a\) and the Dehn twist \(\dt_{\rho}\). The action of $G_d$ on $\pi_6(X_d)$ is given by  
$$f_a \cdot \dt_{\epsilon} = - \dt_{\epsilon}, \ \dt_{\rho} \cdot \dt_{\epsilon} = \dt_{\epsilon}, \ f_a \cdot f_{\eta} = -f_a{_*}(\eta), \ \dt_{\rho} \cdot f_{\eta} =\dt_{\rho}{_*}(\eta), $$
where $f_a{_*}$ and $\dt_{\rho}{_*}$ stand for the induced isomorphism on $\pi_6(X_d)$.

\begin{thm}\label{thm:mainnonspin}
Let $M$ be a simply connected non-spin homology $\cp^3$. 

If $d \ne 0$, then there are short exact sequences
$$1 \to H \to \e(M) \to \mathbb Z/2 \to 1$$
$$1 \to G \to H \to \mathbb Z/2 \to 1.$$
The group $G$ is isomorphic to $\z/4 \times \z/(6,d/2)$,
the action of $\z/2$ on $G$ is given by multiplication by $-1$.
 
 If $d=0$, then there is a short exact sequence
 $$0 \to (\z/4 \times \z/6) \rtimes (\z/2 \times \z/2) \to \e(M) \to \z/2 \to 1.$$
 In the semi-direct product, the actions of $(1,0)$ and $(0,1)$ on $\z/4 \times \z/6$ are both given by multiplication by $-1$.
\end{thm}

The methods used in this paper are primarily drawn from manifold topology. In \S \ref{sec:mcg}, we first recall the results on the Torelli group \(\T(M)\) from \cite{KrSu}, and then extend them to the full mapping class group \(\m(M)\). In \S \ref{sec:forget}, we study the forgetful homomorphism
$\psi_0 \colon \T(M) \to \e_0(M)
$
using the Browder–Novikov–Sullivan–Wall surgery exact sequence. Finally, in \S \ref{sec:main}, we state and prove Theorem \ref{thm:hmtp1} and \ref{thm:hmtp2}. 

\medskip 

\noindent \emph{Acknowledgement.}  The second author is partially supported by NSFC 12471069.

\section{The full mapping class group $\m(M)$}\label{sec:mcg}
\subsection{The Torelli group $\T(M)$}
The Torelli group \(\T(M)\) is an important subgroup of the mapping class group \(\m(M)\). It consists of isotopy classes of diffeomorphisms of \(M\) that induce the identity on \(H_*(M)\). In \cite{KrSu}, the Torelli groups of a more general class of 6-manifolds were computed, namely, simply connected 6-manifolds whose homology groups are isomorphic to
$
H_*(\cp^3 \# r(S^3 \times S^3)).
$
In this subsection, we review this computation in the case of simply connected homology \(\cp^3\)'s. For further details, see \cite[\S 6--8, \S 10]{KrSu}.

Let $M$ be a simply connected homology $\cp^3$, $D \subset M$ be a fixed embedded $6$-disk. Let $\T(M,D)$ be the relative Torelli group, i.e., the group of isotopy classes of diffeomorphisms $f \colon M \to M$ such that $f |_D = \mathrm{id}$ and $f_* = \mathrm{id}$ on $H_*(M)$. It is shown in \cite[Theorem 9.4, \S 10.3]{KrSu} that the homomorphism $\T(M,D) \to \T(M)$ is an isomorphism. Therefore we will identify  these two groups in the sequel.
There is an invariant $S \colon \T(M, D) \to \z^3/ L(M)$ constructed from the mapping torus $M_{f} = M \times [0,1] / (x,0) \sim (f(x),1)$, which is shown to be an isomorphism. 

The definitions of the invariant are slightly different for spin and nonspin manifolds. 
We start with spin manifolds. 
Let $B=\cp^{\infty} \times K(\z,4) \times B\mathrm{Spin}$, $p \colon B \to B\mathrm{Spin} \to  BSO$ be the canonical projection. Then the normal Gauss map $\nu \colon M_{f} \to BSO$ has a lifting $g \colon M_f \to B$ constructed as follows. Fix an orientation of $M$ and a generator $x \in H^2(M)=H^{2}(M_{f})$,  take a map $M_{f} \to \cp^{\infty}$ representing $x$. For $d \ne 0$, let $y=x^{2}/d$ be the generator of $ H^{4}(M)$; for $d=0$, $l \ne 0$, let $y=p_{1}(M)/24l \in H^{4}(M)$; for $d=l=0$, take $y \in H^{4}(M)$ an arbitrary generator. Since $f$ induces the identity on $H_*(M)$, the Wang sequence implies that $H^4(M_f) \to H^4(M)$ is an isomorphism. Let $y' \in H^{4}(M_{f})$ be the preimage of $y$, then $x^{2}=dy'$ and $p_1(M_{f})=4(d+6l)y'$.   Take a map $M_{f} \to K(\z,4)$ representing $y'$. Finally choose the spin structure on the stable normal bundle of $M_{f}$, which restricts to the trivial spin structure on the preferred embedding $D \times S^{1} \subset M_{f}$. This induces a map $M_{f} \to B\mathrm{Spin}$. Putting all these maps together we get $g \colon M_{f} \to B$ which is a lift of the normal Gauss map $M_{f}\to BSO$. Note that all maps involved are unique up to homotopy in the appropriate sense. 

The pair $(M_f, g)$ represents a bordism class  in the $B$-bordism group $\Omega_{7}(B,p)= \omsp_{7}(\cp^{\infty} \times K(\z,4))$. By a  standard computation with the Atiyah-Hirzebruch spectral sequence one shows that this group is zero (\cite[Theorem 6]{Kre18}). Therefore there exists an $8$-dimensional normal $B$-bordism $G \colon W \to B$ such that $\partial (W, G)=(M_{f}, g)$. After taking connected sum with quaternionic  projective planes $\hp^2$'s (or $\hp^2$'s with opposite orientation) we may assume that the signature of $W$ is zero.

Let $\ov x \in H^2(B)$ be the pull-back of the canonical generator of $H^{2}(\cp^{\infty})$, $\ov y \in H^4(B)$ be the pull-back of the canonical generator  of $H^4(K(\z,4))$ and $\frac{p_1}{2} \in H^4(B)$ be the pull-back of the universal spin Pontrjagin class in $H^4(B\mathrm{Spin})$. Let 
$$\alpha = G^*(\ov x^2-d\ov y), \  \beta =G^*(\frac{p_1}{2}+2(d+6l)\ov y) \in H^4(W).$$ 
Then by construction, the restriction of  $\alpha$ and $\beta$ on $\partial W=M_{f}$ are zero, so we have preimages  $\ov{\alpha}$, $\ov{\beta}\in H^4(W, \partial W)$ of $\alpha$ and $\beta$ respectively. Now define characteristic numbers
$$\begin{array}{rcl}
s_1 & = & \langle \ov{\beta}^2, [W,\partial W] \rangle, \\
s_2 & = & \langle \ov{\alpha}^2 + \ov{\alpha}\cup \ov{\beta}, [W,\partial W] \rangle, \\
s_3 & = & \langle \ov{\alpha}\cup \ov{\beta}, [W,\partial W] \rangle. 
\end{array}$$
By comparing $\T(M,D)$ and the mapping class group of $(S^2 \times D^4, \partial)$, and taking consideration of the choices of the $B$-bordism $W$ one has a homomorphism  
$$S = [(s_1/8, s_2/2, s_3)] \colon \T(M, D) \to \z^3 / L(M),$$
where $L(M)$ is the lattice generated by the characteristic numbers of elements in the subgroup of the $8$-dimensional $B$-bordism group $\Omega_8(B)=\omsp_8(\cp^{\infty} \times K(\z,4))$, consisting of manifolds with zero signature. Using the Atiyah-Hirzebruch spectral sequence one can show that this group is isomorphic to $\z^6$ with explicit generators $e_1, \cdots, e_6$ (\cite[Theorem 10]{Kre18}). The  characteristic numbers of these generators are given in the following table, where $k=(d+6l)$,
\begin{equation}\label{eq:table}
\begin{array}{c|cccccc}
& e_1 & e_2 & e_3 & e_4 & e_5 & e_6 \\
\hline
s_{1}/8 & -28 & 0 & 0 & (k^2-k)/2 & k^2 & 0 \\
s_{2}/2& 0 & 12 & 0 & (d^2-d(2k-1))/2 & d^{2}-2dk & k-d \\
s_{3} & 0 & 0 & 2 & d(2k-1) & 4dk & -2k
\end{array}
\end{equation}
 The lattice $L(M) \subset \z^3$ is generated by column vectors.  In \cite[Theorem 6.1]{KrSu} it is shown that $S$ is an isomorphism, and the groups $\T(M)$ are given there. In particular, there is a surjective homomorphism $\z/28 \times \z/12 \times \z/2 \to \T(M)$.
 
\begin{rem}
In the definition of the $S$-invariant, we fix an orientation of $M$ and a generator $x \in H^2(M)$. From the definition it is ready to see that if we reverse the orientation then the invariant changes by a sign, while if we use the generator $-x$ then the invariant remains unchanged. 
\end{rem} 
 
A diffeomorphism $f \colon M \to M$ is \emph{disk supported} if $f= \mathrm{id}$ outside an embedded closed $6$-dimensional disk $D \subset M$. Isotopy classes of disk supported diffeomorphisms form a subgroup $K \subset \T(M)$. The quotient group $\T(M) / K$ will be relevant in the computation of $\e(M)$ (see Proposition \ref{prop:kernel}). Let $\m(D^6, \partial)$ be the mapping classes group of the $6$-disk relative to boundary. Fix an embedding  $D^6 \subset M$, then there is a homomorphism $\m(D^6, \partial) \to \T(M)$, extending a diffeomorphism of $D^6$ by the identity outside the embedded $D^6$ to a diffeomorphism of $M$. By the disk theorem, the image of this isomorphism is exactly the subgroup $K$.  The group $\m(D^6, \partial)$ is isomorphic to the Kervaire-Milnor group $\Theta_7$ of $7$-dimensional homotopy spheres, which is isomorphic to $ \z/28$, with a preferred generator $1$ represented by the boundary of the $E_8$-plumbing. The $S$-invariant of the image of  $1 \in \z/28$ in $\T(M)$ equals $[(1,0,0)]$. Therefore from Table (\ref{eq:table}) we get the quotient group $\T(M) / K$.

\begin{prop}\label{prop:local}
Let $M$ be a simply connected  spin homology $\cp^3$. 
\begin{enumerate}
\item If $d$ is odd, 
$$\T(M) / K \cong \left \{ \begin{array}{ll}
 \z/3 & \mathrm{if} \ d \equiv 0 \pmod 3 \\
0 & \mathrm{if} \ d \not\equiv 0 \pmod 3 \end{array} \right.$$

\item if $d \equiv 2 \pmod 4$, 
$$\T(M) / K \cong \left \{ \begin{array}{ll}
 \z/2 \times \z/3 & \mathrm{if} \ d \equiv 0 \pmod 3 \\
\z/2 &  \mathrm{if} \ d \not\equiv 0 \pmod 3 \end{array} \right.$$

\item if $d \equiv 4 \pmod 8$, $$\T(M) / K \cong \left \{ \begin{array}{ll}
 \z/2 \times \z/2 \times \z/3 & \mathrm{if} \ d \equiv 0 \pmod 3 \\
\z/2 \times \z/2 & \mathrm{if} \ d \not\equiv 0 \pmod 3 \end{array} \right.$$

\item if $d \equiv 0 \pmod 8$, $\T(M)/K$ is isomorphic to 
$$\begin{array}{c|c|c}
 & d \equiv 0 \pmod 3  & d \not\equiv 0 \pmod 3 \\
\hline
l \mathrm{\ even} & \z/2 \times \z/4 \times \z/3  & \z/2 \times  \z/4 \\

\hline
l \mathrm{\ odd} & \z/2 \times \z/2 \times \z/3  &  \z/2 \times \z/2 \\

\end{array}$$
\end{enumerate}
\end{prop}
\begin{proof}
From Table (\ref{eq:table}) we see  that the finite abelian group $\T(M)/K$ has only $2$- and $3$-torsions, with exponents at most $2$ and $1$, respectively. To determined the exponents we only need to take tensor product with $\z/4$ or $\z/3$, respectively. Then an elementary calculation shows this proposition.
\end{proof}

When $M$ is nonspin, we still take $B=\cp^{\infty} \times K(\z,4) \times B\mathrm{Spin}$, the fibration $p$ is 
$$p \colon B=\cp^{\infty} \times K(\z,4) \times B\mathrm{Spin}\to BSO \times BSO \stackrel{\oplus}{\longrightarrow} BSO,$$
where $\cp^{\infty} \to BSO$ is the classifying map of the canonical complex line bundle $\eta$ over $\cp^{\infty}$ and $K(\z,4) \times  B\mathrm{Spin} \to  B\mathrm{Spin} \to BSO$ is the composition of natural projections. Then similar to the spin case, for a diffeomorphism $f$, we have a lifting $g \colon M_{f} \to B$ of the normal Gauss map of $M_f$. The corresponding bordism group $\om_{*}(B,p) $ is the twisted spin cobordism group $ \omsp_{*}(\cp^{\infty} \times K(\z,4);\eta)$. This group vanishes (\cite[Theorem 6]{Kre18}).  Let $\ov x \in H^2(\cp^{\infty})$, $\ov y \in H^4(K(\z,4))$ and $\frac{p_1}{2} \in H^4(B\mathrm{Spin})$ be as in the spin case. Let $G \colon W^8 \to B$ be a coboundary of $g \colon M_{f} \to B$ with signature zero, let $\alpha =G^*(\ov x^2-d \ov y)$, $\beta = G^*(\frac{p_1}{2}+(d+24l)\ov y) \in H^4(W)$, then they have preimages $\ov{\alpha}$, $\ov{\beta}\in H^4(W, \partial W)$ respectively. 
Define characteristic numbers
$$ \begin{array}{rcl}
\sigma_1 & = & \langle \ov \alpha \cup \ov \alpha ,[W,\partial W] \rangle \\
\sigma_2 & = & \langle \ov \alpha \cup \ov \beta ,[W,\partial W] \rangle \\
\sigma_3 & = & \langle \ov \beta \cup \ov \beta ,[W,\partial W] \rangle . \\
\end{array} $$
Let $\lambda=\sigma_1+3\sigma_2+\sigma_3$, $\mu=(6\sigma_1+10\sigma_2+3\sigma_3)/4$, then the $S$-invariant  
$$S=[(\mu, \lambda)] \colon \T(M, D) \to \z^{2}/ L(M)$$
 is an isomorphism, where $L(M)$ is the lattice 
 generated by the column vectors of the following table.
$$\begin{array}{c|ccccc}
\mu & 6\cdot 28 & 0 & -\delta^2-48\delta l +2\delta +432l^2 -36l & -2\delta^2-96\delta l -864l^2& 3(\delta-4l)\\
\lambda & 0 & -4 & -4\delta^2-48\delta l +2\delta +576l^2-48l & -2(4\delta^2 +48\delta l -576l^2) & 6(\delta -4l)
\end{array}$$
The image of  the preferred generator of $\m(D^6, \partial) \cong \z/28$ in $\T(M)$ has $S$-invariant $[(-6,-12)]$.  From this one may compute the quotient group $\T(M) / K$.
\begin{prop}\label{prop:local2}
Let $M$ be a simply connected  nonspin homology $\cp^3$, $\T(M)/K$ is given in the following table. 
$$\begin{array}{c|c|c|c|c}

  d & 12 | d & \ \ 3 | d, 4 \nmid d \ \ & 3 \nmid d, 4 | d & 3 \nmid d, 4 \nmid d  \\
\hline
 \T(M)/K & \mathbb Z/6 \times \z/4 & \mathbb Z/3 \times \z/4 & \mathbb Z/2 \times \z/4 & \z/4 \\

 \end{array}$$ 
\end{prop}

\subsection{Dehn twists}\label{subsec:dehn}
In this subsection we describe generators of the group $\T(M)/K$.  

Let $i \colon S^2 \to M$ represent a generator of $\pi_2(M) \cong \z$, $\varphi \colon S^3 \to S^2$ be the Hopf fiber map. Then $i \circ \varphi \colon S^3 \to M$ represents a generator of $\pi_3(M) \cong \z$. Fix a thickened embedded $3$-sphere $S^3 \times D^{3} \subset M$, such that $S^3 \times \{0\}$ represents this generator of $\pi_3(M)$. For any  $\alpha \in \pi_{3}SO(4)$, choose a representing smooth map $f \colon D^{3} \to SO(4)$, which maps an open neighborhood of $\partial D^{3}$ to $I$. The \emph{Dehn twist} 
$\mathbf{t}_{\alpha} \colon M \to M$ is defined 
as the diffeomorphism which maps the complement of $S^3 \times D^{3}$ by the identity and $(x,y) \in S^3 \times D^{3}$ to $(f(y)x, y)$. 

From the exact sequence of homotopy groups of fiber bundle 
$$SO(3) \to SO(4) \to S^3$$
it can be shown that the group $\pi_3(SO(4))$ is isomorphic to $\z \oplus \z$, with generators represented by $\epsilon \colon S^3 \to SO(3) \to SO(4)$, which is the composite of the projection $S^3 \to SO(3)$ with the inclusion,  and $\rho \colon S^3 \to SO(4)$, which is the inclusion of the unit quaternions $S^3$ to $SO(4)$. The corresponding Dehn twists are denoted by $\dt_{\epsilon}$ and $\dt_{\rho}$, respectively.

\begin{lem}\label{lem}
\begin{enumerate}
\item The group $\T(M) /K$ is generated by $\dt_{\epsilon}$ and $\dt_{\rho}$; the order of $\dt_{\rho}$ is at most $2$.
\item
Let $J \colon \pi_3SO(3) \to \pi_6(S^3)$ be the $J$-homomorphism, $\xi=i \circ \varphi \circ J(\epsilon) \colon S^6 \to M$. 
Then the homotopy equivalence $f_{\xi}=(\mathrm{id} \vee \xi )\circ \mu \colon M \to M \vee S^6 \to M$
is homotopic to $\dt_{\epsilon}$. 

\end{enumerate}
\end{lem}
\begin{proof}
The statements in (1) is shown in \cite[\S 7]{KrSu}. To show (2), consider the natural inclusion map $SO(3) \to G^+(S^3)$, where $G^+(3)$ is the space of base point and orientation preserving homotopy equivalences of $S^3$. The induced homomorphism is the $J$-homomorphism
$$J \colon \pi_3SO(3) \to \pi_3G^+(3) \cong \pi_6(S^3).$$
Therefore regarded as a diffeomorphism of $S^3 \times D^3$, the Dehn twist $\dt_{\epsilon}$ is homotopic to 
$$(\mathrm{id} \vee J(\epsilon) )\circ \mu \colon S^3 \times D^3 \to (S^3 \times D^3) \vee S^6 \to S^3 \times D^3,$$ where the map on $S^6$ is the composition of  $J(\epsilon) \colon S^6 \to S^3 \times \{0\}$ and the inclusion $S^3 \times \{0\} \to  S^3 \times D^3$. As a diffeomorphism of $M$, the Dehn twist $\dt_{\epsilon}$ is homotopic to $(\mathrm{id} \vee \xi )\circ \mu$.
\end{proof}

\subsection{The full mapping class group $\m(M)$}
In this subsection we compute the full mapping class group $\m(M)$, in which orientation reversing diffeomorphims are contained as well.
\begin{thm}\label{thm:full1}
Let $M$ be a simply connected homology $\cp^3$ which is not diffeomorphic to $S^2 \times S^4$ or $S^4 \widetilde{\times} S^2$. Then there is a short exact sequence \begin{equation}\label{seq:1}
1 \to \T(M) \to \m(M) \to \z/2 \to 1,
\end{equation}
where the action of $\z/2$ on $\T(M)$ is by $-1$. This short exact sequence splits when $M$ is spin, hence $\m(M) = \T(M) \rtimes \z/2$.
\end{thm}
\begin{proof}
Denote by $\m^+(M)$ the oriented mapping class group of $M$. Then there is an exact sequence 
$$1 \to \m^+(M) \to \m(M) \to \z/2.$$
It is assumed that $M$ is not diffeomorphic to $S^2 \times S^4$ or $S^4 \widetilde{\times} S^2$. Therefore either $d \ne 0$ or $l \ne 0$. If $d \ne 0$, then clearly $\T(M) = \m^+(M)$. If $d=0$ but $l \ne 0$, then $p_1(M) \ne 0$. In this case since a diffeomorphism  preserves $p_1(M)$, it induces the identity on $H^4(M)$, therefore $\T(M)=\m^+(M)$. Moreover, by the classification theorem (Theorem \ref{thm:class}), there exists an orientation-reversing self-diffeomorphism of $M$, therefore we  have a  short exact sequence
\begin{equation*}
1 \to \T(M) \to \m(M) \to \z/2 \to 1.
\end{equation*}
Since $\T(M)$ is an abelian group, there is an action of $\z/2$ on $\T(M)$, given by conjugation: let $f $ be an orientation-reversing diffeomorphism, $[g] \in \T(M)$, then the action of $1 \in \z/2$ on $[g]$ is $[f^{-1} \circ g \circ f] \in \T(M)$. Observe that there is an orientation-preserving diffeomorphism between the mapping tori
$$M_{f^{-1} \circ g \circ f} \stackrel{\cong}{\longrightarrow} M_{g^{-1}},\ \ \ [(x,t)] \mapsto [(f(x),1-t)]$$
This implies that $f^{-1} \circ g \circ f$ and $g^{-1}$ have the same $S$-invariants, therefore $[f^{-1} \circ g \circ f] = [g^{-1}] \in  \T(M)$. This shows that the action of $\z/2$ on $\T(M)$ is by $-1$.

Next we show that when $M$ is spin the exact sequence is split. If $l$ is even, it's known that there exists a free involution on $M$ (see \cite[Theorem 1.2]{Olber} and the discussion there). By the Lefschetz fixed-point theorem, this free involution is orientation-reversing. This gives a splitting. If $d$ is even, then by \cite[Corollary 1.5]{Olber1}, there exists a conjugation on $M$, which is a special type of involution. Such a conjugation has a $3$-dimensional fixed-point set, hence must be orientation-reversing. This gives a splitting. If both $d$ and $l$ are odd,  then from the tables in \cite[Theorem 6.1]{KrSu} one sees that the Torelli group $\T(M)$ has odd order. Therefore $H^{2}(\z/2; \T(M))=0$. This implies that the sequence (\ref{seq:1}) splits. 
\end{proof}

\begin{thm}\label{thm:full2}
Let $M= S^2 \times S^4$ or $S^4 \widetilde{\times} S^2$, then 
$$\m(M) = \T(M) \rtimes (\z/2 \times \z/2)$$
where the actions of the standard generators of $\z/2 \times \z/2$ on $\T(M)$ are by $-1$.
\end{thm}
\begin{proof}
If $M=S^2 \times S^4$, then clearly we have a short exact sequence
$$1 \to \T(M) \to \m(M) \to \z/2 \times \z/2 \to 1$$
where the homomorphism $\m(M) \to \z/2 \times \z/2$ is given by the actions on $H_2(M)$ and $H_4(M)$. Apparently the elements $(1,0)$ and $(0,1) \in \z/2 \times \z/2$ have liftings $(-\mathrm{id}, \mathrm{id}) \colon S^{2} \times S^{4} \to S^{2} \times S^{4}$ and $(\mathrm{id}, -\mathrm{id}) \colon S^{2} \times S^{4} \to S^{2} \times S^{4}$ respectively. Since these are orientation reversing diffeomorphisms, as discussed in the proof of last theorem,  their actions on $\T(M)$ are by $-1$.

When $M= S^4 \wt S^2$, the non-trivial linear $S^4$-bundle over $S^2$, we construct two smooth involutions as follows. The first involution $f_1$ is the fiberwise antipodal map, then $f_{1}$ induces the identity on $H_2$ and $-\mathrm{id}$ on $H_4$. To construct the second involution $f_2$, we use the decomposition $S^4 \wt S^2 = (S^4 \times D^2)_1 \cup_{\alpha} (S^4 \times D^2)_2$, with $\alpha$ a map $S^1 \to SO(5)$ representing the nontrivial element in $\pi_1SO(5) \cong \z/2$. Consider the space
$$SO(5)^{\z/2} = \{A \in SO(5) \ | \ A^t=A\}.$$
Such matrices have eigenvalues $\pm 1$. It's not hard to see that this space is a union of Grassmannians of even dimensional linear subspaces in $\mathbb R^{5}$. Especially we have an inclusion $G_4(\mathbb R^5) \cong \mathbb R \mathrm P^4 \subset SO(5)^{\z/2} \subset SO(5)$, inducing isomorphism on $\pi_1$. We may choose $\alpha$ to be a loop in $SO(5)^{\z/2}$. Then we may define an involution $f_2 \colon S^4 \wt S^2 \to S^4 \wt S^2$, by interchanging the two copies $S^4 \times D^2$: $(x,y)_1 \mapsto (x,y)_2$, $(x,y)_2 \mapsto (x,y)_1$. Since for any $y \in S^{1}$, $\alpha(y)^{2}=I$, this is a well-defined diffeomorphism. Clearly $f_2$ induces $-\mathrm{id}$ on $H_2$ and the identity on $H_4$, and $f_1$ and $f_2$ commute. Therefore we have $\m(S^4 \wt S^2)=\T(S^4 \wt S^2) \rtimes (\z/2 \times \z/2)$, with the actions of $(1,0)$ and $(0,1)$ are by $-1$.
\end{proof}

\section{The forgetful homomorphism}\label{sec:forget}
The full mapping class group $\m(M)$ and the group of self-equivalences $\e(M)$ are related by the forgetful homomorphism $\psi \colon \m(M) \to \e(M)$. Let $\e_0(M)$ be the subgroup of $\e(M)$ consisting of homotopy equivalences which act trivially on $H_*(M)$, then $\psi$ restricts to a homomorphism $\psi_0 \colon \T(M) \to \e_0(M)$. In this section we analyze this homomorphism by surgery theory.

\begin{prop}\label{prop:hmtp}
Let $M$ be a simply connected homology $\cp^3$, $\psi_0 \colon \T(M) \to \e_0(M)$ be the forgetful homomorphism. When $d$ is odd, $\psi_0$ is surjective; when $d$ is even, the image of $\psi_0$ is an index $2$ subgroup of $\e_0(M)$. 
\end{prop}
\begin{proof}
Denote by $\mathscr S(M)$ the smooth structure set of $M$. Elements of $\mathscr S (M)$ are equivalence classes of homotopy equivalences $f \colon X \to M$ from a closed smooth manifold $X$ to $M$. Two pairs $(X_1, f_1)$ and $(X_2, f_2)$ are equivalent if there exists a diffeomorphism $\varphi \colon X_1 \to X_2$ such that $f_1$ is homotopic to $f_2 \circ \varphi$.  There is a map $\alpha \colon \e_0(M) \to \mathscr S(M)$ sending a homotopy equivalence $f \colon M \to M$ to $[M,f] \in \mathscr S(M)$. By definition, $\alpha(f) = \alpha(g)$ if and only if $f$ and $g$ are in the same coset of $\mathrm{Im}\, \psi_0$. Therefore we have an injective map between sets
$$\ov{\alpha} \colon \e_0(M)/\mathrm{Im}\, \psi_0 \to \mathscr S(M)$$
 
 The image of this map can be analyzed via the surgery exact sequence (\cite[Theorem 14.8]{KL05})
$$ L_7(\z) \to \mathscr S(M) \stackrel{\eta}{\rightarrow} \mathcal N (M) \stackrel{\theta}{\rightarrow} L_6(\z).$$
The surgery obstruction group $L_7(\z)$ is zero, hence $\eta$ is injective. Elements in $\mathcal N(M)$ are represented by degree $1$ normal maps
$$\xymatrix{
\nu N \ar[r]^{\bar f} \ar[d] & E \ar[d] \\
N \ar[r]^f & M.
}$$
More precisely, the set $\mathcal N(M)$ can be described as 
$$\mathcal N(M) =  \oplus \Omega_6(M;E)_1 / \mathrm{Aut}(E),$$
where the direct sum takes over all isomorphism classes of oriented stable vector bundles $E$ over $M$, for the details see \cite[\S 14.5]{KL05}.

We first analyze the group $\mathrm{Aut}(E)$, which is the group of homotopy classes of bundle isomorphisms of $E$. Fix a $k$-dimensional representative $\xi$ of $E$ with $k >>0$, let $P \to M$ be the $SO(k)$-principal bundle associated to $\xi$. Then the group $\mathrm{Aut}(\xi)$ is identified with $\pi_0(\mathcal G(P))$, where $\mathcal G(P)$ is the gauge group of the principal bundle $P$. It is known that the gauge group $\mathcal G(P)$ is homotopy equivalent to  the based loop space $\Omega \mathrm{Map}(M, BSO(k))$ of the space of continuous maps from $M$ to $BSO(k)$, with base point a classifying map $\varphi \colon M \to BSO(k)$ of $\xi$ (see \cite{AB83}). Therefore $\mathrm{Aut}(\xi)$ is isomorphic to $\pi_1(\mathrm{Map}(M, BSO(k)), \varphi)$. From the fiber sequence
$$\Omega BSO(k) \to \mathrm{Map}_*(M, BSO(k)) \to \mathrm{Map}(M, BSO(k)) \stackrel{\mathrm{ev}}{\longrightarrow} BSO(k)$$
where $\mathrm{Map}_*$ denotes the space of based continuous maps, and $\mathrm{ev}$ is the evaluation map at the base point $x_0 \in M$, 
one has an exact sequence
\begin{equation}\label{eq:map}
 \pi_1(\mathrm{Map}_*(M, BSO(k)), \varphi) \to \pi_1(\mathrm{Map}(M, BSO(k)), \varphi) \to \pi_0SO(k).
\end{equation}
By the adjunction formula, the space $\mathrm{Map}_*(S^1, \mathrm{Map}_*(M, BSO(k))_{\varphi})$ is homeomorphic to the space of continuous maps from $S^1 \times M/ S^1 \times \{x_0\}$ to $BSO(k)$ whose restriction on $\{s_0\} \times M$ is $\varphi$, where $s_0 \in S^1$ is the base point. Notice that all the groups  
$$H^i(S^1\times M / S^1 \times \{x_0\}, \{s_0\} \times M; \pi_iBSO(k)) \cong H^i(S^1 \wedge M; \pi_iBSO(k))$$ vanish. Elementary obstruction theory shows that all such maps are homotopic relative to $\{s_0\} \times M$. Therefore the exact sequence (\ref{eq:map}) implies $\mathrm{Aut}(\xi)=0$. By passing to the limit one has $\mathrm{Aut}(E)=0$.  Hence $\mathcal N(M) = \oplus \Omega_6(M;E)_1$.

Next we analyze the map $\eta$.
The element $\eta([M,f]) \in \mathcal N(M)$ is represented by the degree $1$ normal map $(\bar f ,f)$
$$\xymatrix{
\nu M \ar[r]^{\bar f} \ar[d] & g^*(\nu M) \ar[d] \\
M \ar[r]^f & M 
}$$ 
where $g$ is a homotopy inverse of $f$. 
By the minimal cell structure theorem, $M$ is homotopy equivalent to $P_{d} \cup e^{6}$, where $P_{d} = S^{2} \cup_{d\eta} e^{4}$, $\eta \colon S^{3} \to S^{2}$ the Hopf map. The Puppe sequence of the cofibration sequence
$S^{5} \to P_{d} \to M \to S^{6}$,
together with the fact $\pi_5(BO)=\pi_6(BO)=0$ implies that the restriction map $[M, BO] \to [P_d,BO]$ is an isomorphism. A further analysis of the Puppe sequence of the cofibration sequence $S^3 \to S^2 \to P_d \to S^4$ shows that a stable vector bundle $\xi$ over $P_d$ is determined by $w_2(\xi)$ and $p_1(\xi)$. Since $g$ acts trivially on homology, $g^*(\nu M)$ is isomorphic to $\nu M$. Therefore $\eta([M,f])$ is in the subset $\Omega_6(M;\nu M)_1$ of $\mathcal N(M)$.

The bordism group $\Omega_6(M;\nu M)$ is identified by the Thom-Pontrjagin construction with the homotopy group $\pi_6(\mathrm{Th}(\nu M))$ of the Thom spectrum $\mathrm{Th}(\nu M)$, whose $n$-th space is the Thom space of the normal bundle of $M \subset \mathbb R^{6+n}$.  We compute the homotopy group of the spectrum by the Atiyah-Hirzebruch spectral sequence. The $E_2$-terms $E^2_{p,q}$ are identified via the Thom isomorphism with $H_p(M;\pi_q^S)$. The differential $d_2 \colon E^2_{6,1} \to E^2_{4,2}$ is dual to 
$$Sq^2+w_2(M)\cup - \colon H^4(M;\z/2) \to H^6(M;\z/2);$$ 
$d_2 \colon E^2_{6,0} \to E^2_{4,1}$ is the composition of the mod $2$ reduction and the dual of $Sq^2 +w_2(M)\cup -\colon H^4(M;\z/2) \to H^6(M;\z/2)$ (\cite[Lemma]{Te93}). Using the same idea in the proof of Lemma in \cite{Te93}, it can be shown that the differential $d_2 \colon E^2_{4,2} \to E^2_{2,3}$ is the composition of the dual of $Sq^2 +w_2(M)\cup - \colon H^2(M;\z/2) \to H^4(M;\z/2)$ and the homomorphism $j_* \colon H_2(M,\z/2) \to H_2(M;\z/24)$ induced by the inclusion of coefficients. (The differential $d_2$ is a stable homology operation. Such operations are  classified by 
$$[H\z/2, \Sigma^2 H\z/24] \cong H^2(H\z/2; \z/24) \cong \z/2 \{ j_* \circ (Sq^2)^*\},$$
and the non-triviality of $d_2$ can be detected by $\mathbb C\mathrm P^2$, whose stable homotopy groups can be calculated using the cofiber sequence $S^2 \to \mathbb C \mathrm P^2 \to S^4$.) A straightforward calculation shows that $d_2 \colon E^2_{6,1} \to E^2_{4,2}$ and $d_2 \colon E^2_{6,0} \to E^2_{4,1}$ are trivial; $d_2 \colon E^2_{4,2} \to E^2_{2,3}$ is nontrivial when $d$ is odd, and trivial when $d$ is even. 

Therefore, when $d$ is odd, we get a short exact sequence
$$0 \to \Omega_6^{\mathrm{fr}} \to \Omega_6(M;\nu M) \to H_6(M) \to 0.$$
The surgery obstruction map $\theta \colon \Omega_6^{\mathrm{fr}} \to L_6(\z)$ is an isomorphism. 
The normal map $(\bar f, f)$ has vanishing surgery obstruction and  the degree of $f$ is $1$. Therefore the image of $\eta \circ \overline{\alpha}$ consists of one element. Since $\eta$ is injective, this shows that $\e_0(M) / \mathrm{Im} \, \psi_0$ consists of one element. When $d$ is even, there are short exact sequences
$$0 \to G \to \Omega_6(M;\nu M) \to H_6(M) \to 0$$
$$0 \to \Omega_6^{\mathrm{fr}} \to G \to \z/2 \to 0.$$
Therefore there are two bordism classes of degree one normal maps with trivial surgery obstruction. Such an element is represented by 
$$\xymatrix{
\nu N \ar[r]^{\bar g} \ar[d] & \nu M \ar[d] \\
N \ar[r]^g & M 
}$$ 
with $g$ a homotopy equivalence. By the classification of homology $\cp^3$ (Theorem \ref{thm:class}), $N$ is diffeomorphism to $M$. Therefore, the image of $\eta \circ \overline{\alpha}$ consists of two elements. This shows that $\mathrm{Im}\, \psi_0$ is a subgroup of index $2$ in $\e_0(M)$.
\end{proof}

\begin{prop}\label{prop:kernel}
The kernel of $\psi_0$ is the subgroup of $\T(M)$ consisting of disk-supported diffeomorphisms.
\end{prop}
\begin{proof}
Clearly any disk-supported diffeomorphism is in the kernel of $\psi_0$ by the Alexander's trick. On the other hand, given a diffeomorphism $f$ which is homotopic to the identity, let $h \colon M \times [0,1] \to M \times [0,1]$ be a homotopy, with $h_0=\mathrm{id}$ and $h_1=f$, then $[M \times [0,1] ,h]$ can be viewed as an element in the relative structure set $\mathscr S_{\partial}(M \times [0,1])$. Consider the surgery exact sequence
$$L_8(\z) \to \mathscr S_{\partial}(M \times [0,1]) \to [\Sigma M, G/O] \to L_7(\z).$$
The homomorphism $L_8(\z) \to \mathscr S_{\partial}(M \times [0,1])$ is given by taking  connected sum of some $7$-dimensional homotopy sphere $\Sigma$ with the base point $[M \times [0,1], \mathrm{id}]$ in $\mathscr S_{\partial}(M \times [0,1])$.
By obstruction theory it's easy to see that $[\Sigma M, G/O]=0$ (since $\pi_i(G/O)=0$ for $i=1,3,5,7$).
Therefore there is a diffeomorphism $g \colon M \times [0,1] \to (M \times [0,1])\# \Sigma$ such that $g|_{M \times 0}=\mathrm{id}$, $g|_{M \times 1}=f$. Composed with the diffeomorphism $\varphi \colon (M \times [0,1]) \# \Sigma \to M \times [0,1]$, where $\varphi|_{M \times 0}$ is the disk supported diffeomorphism corresponding to $\Sigma$ and $\varphi|_{M \times 1}=\mathrm{id}$, we have a diffeomorphism $k \colon M \times [0,1] \to M \times [0,1]$ such that $k|_{M \times 0}$ is the disk supported diffeomorphism corresponding to $\Sigma$, and $k|_{M \times 1}=f$. This shows that $f$ is pseudo-isotopic to a disk supported diffeomorphism. By the pseudo-isotopy theorem of Cerf \cite{Cerf}, $f$ is isotopic to a disk supported diffeomorphism. 

For a more general discussion of the homomorphism $\T(M) \to \e_0(M)$ in  higher dimensions see \cite[Proposition II.9]{Brum}.
\end{proof}
\section{Determination of $\e(M)$}\label{sec:main}
Let $\psi \colon \m(M) \to \e(M)$ be the forgetful homomorphism, $\e_0(M)$ be the group of homotopy equivalences which acts trivially on homology.
\begin{thm}\label{thm:hmtp1}
Let $M$ be a simply connected spin homology $\cp^3$. 
\begin{enumerate}
\item If $d$ is odd, then $\psi$ is surjective,  
$$\e(M) = \e_0(M)\rtimes \z/2, \ \ \e_0(M) \cong \z/(3,d),$$
the action of $\z/2$ on $\e_0(M)$ is by $-1$.

\item If $d$ is even and $d \ne 0$,  then there is a short exact sequence
$$1 \to \mathrm{Im} \, \psi \to \e(M) \to \z/2 \to 1,$$ 
 where $\mathrm{Im} \, \psi$ is isomorphic to $G \rtimes \z/2$, and
 $$G \cong \left \{ \begin{array}{ll}
\z/2 \times \z/(6, d/2) & \mathrm{if} \ d \not\equiv 0 \pmod 8 \ \mathrm{and} \ l \not\equiv 0 \pmod 2 \\
\z/2 \times \z/(12, d/2) & \mathrm{otherwise}
\end{array} \right.$$
The action of $\z/2$ on $G$ is by $-1$, and $G \subset \e_0(M)$ is a subgroup of index $2$. 
\item If $d=l= 0$, i.e., $M = S^2 \times S^4$, then there is a short exact sequence
$$1 \to \mathrm{Im} \, \psi \to \e(M) \to \z/2 \to 1,$$ 
where $\mathrm{Im} \, \psi$ is isomorphic to $G \rtimes (\z/2 \times \z/2)$, $G \cong \z/2 \times \z/12$, the actions of $(1,0)$ and $(0,1)$ on $G$ are by $-1$. $G \subset \e_0(M)$ is a subgroup of index $2$.

\item If $d=0$, $l$ is even and non-zero, then $\mathrm{Im} \, \psi$ is a subgroup of index $4$, isomorphic to $G \rtimes \z/2$, $G \cong \z/2 \times \z/12$, the action of $\z/2$ on $G$ is by $-1$. $G \subset \e_0(M)$ is a subgroup of index $2$.

\item If $d=0$ and $l$ is odd, then $\mathrm{Im} \, \psi$ is a subgroup of index $2$ or $4$, isomorphic to $(\z/2 \times \z/6) \rtimes \z/2$,  the action of $\z/2$  is by $-1$.
\end{enumerate}
\end{thm}

\begin{thm}\label{thm:hmtp2}
Let $M$ be a simply connected nonspin homology $\cp^3$. 
\begin{enumerate}
\item If $d \ne 0$, then there are short exact sequences
$$1 \to \mathrm{Im} \, \psi \to \e(M) \to \z/2 \to 1$$
$$1 \to G \to \mathrm{Im} \, \psi \to \z/2 \to 1,$$
where $G$ is isomorphic to 
$$\begin{array}{c|c|c|c|c}

  & 12 | d & \ \ 3 | d, 4 \nmid d \ \ & 3 \nmid d, 4 | d & \mathrm{otherwise}  \\
\hline
\ \ \ G  \ \ \ & \mathbb Z/6 \times \z/4 & \mathbb Z/3 \times \z/4 & \mathbb Z/2 \times \z/4  & \z/4  \\

 \end{array}$$ 
The action of $\z/2$ on $G$ is by $-1$. $G \subset \e_0(M)$ is a subgroup of index $2$.

\item If $d=l=0$, i.e. $M = S^4 \widetilde{\times} S^2$,  then there is a short exact sequence
$$1 \to \mathrm{Im} \, \psi \to \e(M) \to \z/2 \to 1,$$
where $\mathrm{Im} \, \psi $ is isomorphic to $G \rtimes (\z/2 \times \z/2)$, $G \cong \z/6 \times \z/4$,  the actions of $(1,0)$ and $(0,1)$ on $G$ are by $-1$. $G \subset \e_0(M)$ is a subgroup of index $2$.

\item If $d=0$ and $l \ne 0$, then $\mathrm{Im} \, \psi $ is a subgroup of index $4$, and there is a short exact sequence
$$ 0 \to G \to \mathrm{Im} \, \psi  \to \z/2 \to 0,$$
where $G \cong \z/6 \times \z/4$, the action of $\z/2$ on $G$ is by $-1$. $G \subset \e_0(M)$ is a subgroup of index $2$
\end{enumerate}
\end{thm}

\begin{proof}[Proof of Theorem \ref{thm:hmtp1} and \ref{thm:hmtp2} ]
First assume $d \ne 0$. Then by Theorem \ref{thm:full1} we have the following commutative diagram, where rows and columns are exact sequences
$$\xymatrix{
& \ker \, \psi_0 \ar[d] \ar[r]^= & \ker \, \psi \ar[d] & \\
0  \ar[r] &  \T(M) \ar[d]^{\psi_0} \ar[r] & \m(M) \ar[d]^{\psi} \ar[r] & \z/2 \ar[d]^= \ar[r] & 1\\
0 \ar[r] & \e_0(M)  \ar[r] & \e(M)  \ar[r] & \z/2 \ar[r] & 1 }$$

When $M$ is spin, by Theorem \ref{thm:full1} the upper exact sequence splits. 
If $d$ is odd, by Proposition \ref{prop:hmtp} $\psi_0$ is surjective. This implies that $\psi$ is surjective, and $\e(M) = \e_0(M) \rtimes \z/2$, with $\e_0(M) = \mathrm{Im} \, \psi_0$. If $d$ is even and $d \ne 0$, by Proposition \ref{prop:hmtp} $\mathrm{Im} \, \psi_0$ is a subgroup of index $2$ in $\e_0(M)$. Hence $\mathrm{Im} \, \psi$ is a subgroup of index $2$ in $\e(M)$, and isomorphic to $\mathrm{Im} \, \psi_0 \rtimes \z/2$.  By Proposition \ref{prop:kernel} $\mathrm{Im } \, \psi_0$ is isomorphic to $\T(M)/K$, whose isomorphism class is given in  Proposition \ref{prop:local}. 

When $M$ is non-spin, $d$ is even. By Proposition \ref{prop:hmtp} $\mathrm{Im} \, \psi_0$ is a subgroup of index $2$ in $\e_0(M)$. Hence $\mathrm{Im} \, \psi$ is a subgroup of index $2$ in $\e(M)$. There is a short exact sequence
$$1 \to \mathrm{Im} \, \psi_0 \to \mathrm{Im} \, \psi \to \z/2 \to 1$$
with $ \mathrm{Im} \, \psi_0 \cong \T(M) / K$ given in Proposition \ref{prop:local2}. 

Now consider the cases with $d=0$. 

First assume $l=0$. Then  $M=S^2 \times S^4$ or $S^4 \wt S^2$. By Theorem \ref{thm:full2} we have a similar commutative diagram, 
$$\xymatrix{
& \ker \, \psi_0 \ar[d] \ar[r]^= & \ker \, \psi \ar[d] & \\
0  \ar[r] &  \T(M) \ar[d]^{\psi_0} \ar[r] & \m(M) \ar[d]^{\psi} \ar[r] & \z/2 \times \z/2 \ar[d]^= \ar[r] & 1\\
0 \ar[r] & \e_0(M)  \ar[r] & \e(M)  \ar[r] & \z/2 \times \z/2 \ar[r] & 1 }$$
where the homomorphisms from $\m(M)$ or $\e(M)$ to $\z/2 \times \z/2$ are induced by the actions on $H_2(M)$ and $H_4(M)$, and the two rows are split short exact sequences. By Proposition \ref{prop:hmtp} $\mathrm{Im} \, \psi_0$ is a subgroup of index $2$ in $\e_0(M)$,  hence $\mathrm{Im} \, \psi$ is a subgroup of index $2$ in $\e(M)$, isomorphic to $\mathrm{Im} \, \psi_0 \rtimes (\z/2 \times \z/2)$, with  $\mathrm{Im } \, \psi_0 \cong  \z/2 \times \z/12$ if $M = S^2 \times S^4$, or  $\mathrm{Im } \, \psi_0 \cong \z/6 \times \z/4$ if $M= S^4 \wt S^2$.

If $M$ is spin with $l$ even and $l \ne 0$, or if $M$ is nonspin with $l \ne 0$, then we have the following  commutative diagram
\begin{equation}\label{diag2}
\xymatrix{
& \ker \psi_0 \ar[d] \ar[r] & \ker \psi \ar[d] & \\
0  \ar[r] &  \T(M) \ar[d]^{\psi_0} \ar[r] & \m(M) \ar[d]^{\psi} \ar[r] & \z/2 \ar[d] \ar[r] & 1\\
0 \ar[r] & \e_0(M)  \ar[r] & \e(M)  \ar[r] & \z/2 \times \z/2 \ar[r] & 1 }
\end{equation}
where the first row is by Theorem \ref{thm:full1}, the second row by the fact that $M$ is homotopy equivalent to $S^2 \times S^4$ or $S^4 \wt S^2$, and the homomorphism $\z/2 \to \z/2 \times \z/2$ maps $1$ to $(1,0)$. This shows that 
 $\mathrm{Im} \, \psi$ has index $4$ in $\e(M)$. When $M$ is spin, the upper exact sequence splits (Theorem \ref{thm:full1}), therefore $\mathrm{Im} \, \psi$ is isomorphic to $\mathrm{Im} \, \psi_0 \rtimes \z/2$, with $\mathrm{Im} \, \psi_0 \cong \z/2 \times \z/12$.  When $M$ is nonspin, there is a short exact sequence
 $$1 \to \mathrm{Im} \, \psi_0 \to \mathrm{Im} \, \psi \to \z/2 \to 1$$
 with $\mathrm{Im} \, \psi_0 \cong \z/6$. 

If $M$ is spin and $l$ is odd, it is not known to the authors wether the image of the homomorphism $\e(M) \to \z/2 \times \z/2$ in diagram (\ref{diag2}) is $\z/2$ or $\z/2 \times \z/2$. Therefore  $\mathrm{Im} \, \psi$ is a subgroup of index $2$ or $4$ in $\e(M)$.
\end{proof}

\begin{proof}[Proof of Proposition \ref{prop}]
The first statement is shown in Lemma \ref{lem}. The second statement is shown by examining the codimension $4$ splitting invariant. It suffices to consider the case $M=X_d$, the $S^2$-bundle over $S^4$ with characteristic map $d \in \pi_3SO(3)$. Let $f$ be a smooth map, transverse to a fiber $S^2 \subset X_d$, and homotopic to $f_{\eta}$. Then the preimage $F=f^{-1}(S^2)$ is a framed surface in $X_d$. By inspecting the construction of $f_{\eta}$ one sees that $F$ is a disjoint union of the fiber $S^2$ and $g^{-1}(x_0)$, where $g \colon S^6 \to S^4$ is a smooth map representing the non-trivial element of $\pi_6(S^4)$ and $x_0 \in S^4$ is a regular value of $g$. By the Thom-Pontrjagin construction, the set of concordance classes of framed surfaces in $X_d$ is isomorphic to the homotopy set $[X_d, S^4]$. Therefore in order to show that $f$ is not homotopic to a diffeomorphism, it suffices to show that the map $\pi_6(S^4) \to [X_d, S^4]$ induced by the collapsing map $X_d \to S^6$ is injective.  
From the cofiber sequence
$$S^5 \stackrel{b}{\rightarrow} P_d \to X_d \to S^6 \stackrel{\Sigma b}{\longrightarrow} \Sigma P_d$$
one has an exact sequence $[\Sigma P_d, S^4] \to \pi_6(S^4) \to [X_d, S^4]$. When $d$ is odd, $\Sigma P_d$ is homotopy equivalent to $\Sigma \mathbb C \mathrm P^2$ and it is easy to show $[\Sigma \cp^2, S^4]=0$. When $d$ is even, it is known that  $\Sigma b$ is $d$-multiple of a certain element in $\pi_6(\Sigma P_d)$ (\cite[Lemma 4.5]{Yama}). Since $\pi_6(S^4) \cong \z/2$, the map  $[\Sigma P_d, S^4] \to \pi_6(S^4)$ must be trivial. 

\end{proof}

\begin{proof}[Proof of Theorem \ref{thm:mainspin} and \ref{thm:mainnonspin}]
The statements are rephrasing of Theorem \ref{thm:hmtp1} and \ref{thm:hmtp2}. It remains to show that when $d$ and $l$ are even there is a splitting. Since $f_{\eta}$ represents an element of order $2$ in $\e(M)$, and by Proposition \ref{prop} $f_{\eta}$ is not homotopic to a diffeomorphism, it provides a splitting of the exact sequence 
$$1 \to \m(M) \to \e(M) \to \z/2 \to 1.$$

\end{proof}

\begin{proof}[Proof of Corollary \ref{cor}]
By comparing the cardinalities of the group $G$ in Theorem \ref{thm:mainspin} and the group $\pi_6(X_d)$ in formula (\ref{eq:pi6}), one sees that the map $\pi_6(X_d) \to \e(X_d)$ is injective. Therefore exact sequence (\ref{Yama}) is extended to a short exact sequence
$$0 \to \pi_6(X_d) \to \e(X_d) \to G_d \to 1.$$
Note that homotopy equivalences constructed from $\pi_6(X_d)$ are all orientation-preserving. 
 When $d$ is odd, $G_d \cong \z/2$, a splitting is given by the fiberwise antipodal map $f_a$, which is orientation-reversing.  When $d$ is even, $G_d \cong \z/2 \times \z/2$, a preimage of $(1,0) \in \z/2 \times \z/2$ in $\e(X_d)$ is given by $f_a$. Since in this case $\mathrm{Im} \, \psi$ is isomorphic to $\z/2 \times \z/(12,d/2)$, the order of the Dehn twist $\dt_{\rho}$ must be $2$ (Lemma \ref{lem}), and is not in the image of $\pi_6(X_d)$. It is a preimage of $(0,1) \in \z/2 \times \z/2$.
 \end{proof}
 
\noindent \textbf{Tool Use.} In this paper no instrument, software or AI was use.


\end{document}